\documentclass{article}

\usepackage{amssymb, amsthm, amsmath}
\usepackage{xy} \input xy \xyoption{all}
\usepackage[titlenumbered,algoruled]{algorithm2e}

\theoremstyle{plain}
\newtheorem{theorem}{Theorem}

\theoremstyle{definition}

\numberwithin{equation}{section} 

\newcommand{\footnoteremember}[2]{\footnote{#2}\newcounter{#1}\setcounter{#1}{\value{footnote}}}
\newcommand{\footnoterecall}[1]{{\footnotemark[\value{#1}]}}

\newcommand{\pj}{\mathbb{P}}

\DeclareMathOperator{\PGL}{PGL}
\DeclareMathOperator{\Resultant}{Resultant}

\begin{document}

\title{Deciding trigonality of algebraic curves (extended abstract)}

\author{Josef Schicho and David Sevilla\footnote{Partially supported by Research Project MTM2007-67088 of the Spanish Ministry of Education and Science.}\\RICAM, Austrian Academy of Sciences}

\date{}

\maketitle

Let $C$ be an algebraic curve of genus $g\geq3$. Let us assume that $C$ is not hyperelliptic, so that it is isomorphic to its image by the canonical map $\varphi:C\rightarrow\pj^{g-1}$. Enriques proved in \cite{Enriques1919} that $\varphi(C)$ is the intersection of the quadrics that contain it, except when $C$ is trigonal (that is, it has a $g^1_3$) or $C$ is isomorphic to a plane quintic ($g=6$). The proof was completed by Babbage \cite{Babbage1939}, and later Petri proved \cite{Petri1923} that in those two cases the ideal is generated by the quadrics and cubics that contain the canonical curve. In this context, we present an implementation in Magma of a method to decide whether a given algebraic curve is trigonal, and in the affirmative case to compute a map $C\stackrel{3:1}{\rightarrow}\pj^1$ whose fibers cut out a $g^1_3$. Our algorithm is part of a larger effort to determine whether a given algebraic curve admits a radical parametrization.

\section{Classical results on trigonality}

The following theorem \cite[p. 535]{GriffithsHarris1978} classifies canonical curves according to the intersection of the quadric hypersurfaces that contain them.

\begin{theorem}
For $C\subset\pj^n$ any canonical curve, either
\begin{enumerate}
 \item $C$ is entirely cut out by quadric hypersurfaces; or
 \item $C$ is trigonal, in which case the intersection of all quadrics containing $C$ is the rational normal scroll swept out by the trichords of $C$; or
 \item $C$ is a plane quintic, in which case the intersection of the quadrics containing $C$ is the Veronese surface in $\pj^5$, swept out by the conic curves through five coplanar points of $C$.
\end{enumerate}
\end{theorem}

We can use this to sketch an algorithm to detect trigonality.

\begin{algorithm}
  \caption{Sketch of algorithm to detect trigonality}\label{ALG: sketch detect trigonality}
  \BlankLine
  \KwIn{a non-hyperelliptic curve $C$ of genus $g\geq3$}
  \KwOut{\texttt{true} if $C$ is trigonal, \texttt{false} otherwise}
  \BlankLine
  Compute the canonical map $\varphi:C\rightarrow\pj^{g-1}$ and its image $\varphi(C)$\;
  Compute the intersection $D$ of all the quadrics that contain $\varphi(C)$\;
  \eIf{$D=C$}{
    \Return false\;
  }{
    Determine which type of surface is $D$\;
    \uIf{$D=\pj^2$}{
      \Return true \hspace{5ex}\tcp*[h]{$g=3$}\;
    }
    \uElseIf{$D$ is a rational normal scroll}{
      \Return true\;
    }
    \Else{\Return false \hspace{5ex}\tcp*[h]{Veronese} \;}
  }
\end{algorithm}

\pagebreak

For the computation of the canonical map and the computation of the space of forms of fixed degree containing the image of a polynomial curve, there exist efficient implementations in Magma.

For the identification of the surface $D$, we use the Lie algebra method,
which has been introduced in \cite{Schicho2006} (see also \cite{Schicho2008}) for parametrizing certain classes of Del Pezzo surfaces.

\section{The Lie algebra method}

The Lie algebra of a projective variety is an algebraic invariant which is relatively easy to calculate (it is often cheaper than a Gr\"obner basis of the defining ideal, if only generators are given).

Let $X\subset\pj^n$ be an embedded projective variety. Let $\PGL_{n+1}(X)$ be the group of all projective transformations that map $X$ to itself (this is always an algebraic group). The Lie algebra $L(X)$ of $X$ is defined as the tangent space of $\PGL_{n+1}(X)$ at the identity, together with its natural Lie product. It is a subalgebra of $\mathfrak{sl}_{n+1}$, the Lie algebra of $\pj^n$.

For varieties of general type (in particular curves of genus at least 2), the group $\PGL_{n+1}(X)$ is finite and therefore the Lie algebra is zero. On the other hand, the Veronese surface and the rational scrolls have a Lie algebra of positive dimension. This allows us to reduce the recognition problem for Veronese surfaces/rational scrolls to the computations with Lie algebras and their representations.

If $S$ is a rational normal scroll which is not isomorphic to $\pj^1\times\pj^1$, then $L(S)$ has a Levi subalgebra isomorphic to $\mathfrak{sl}_2$. By decomposing the Lie module given by the representation $\mathfrak{sl}_2\hookrightarrow\mathfrak{gl}_{n+1}$, we can construct a $2\times(n-1)$ matrix $A$, such that the $2\times 2$ minors of $A$ generate the ideal of $S$. The ratio of the two entries of any column defines a map $\rho:S\to\pj^1$ whose fibers are lines. Then the map $C\to\pj^1$ is constructed as
\begin{equation}\label{EQ: isomorphism to P1}
C\stackrel{\varphi}{\rightarrow}\varphi(C)\stackrel{i}{\hookrightarrow}S\stackrel{\rho}{\rightarrow}\pj^1.
\end{equation}

If $S$ is a Veronese surface, then we can construct by similar methods an isomorphism of $C$ with a planar quintic.

\section{Trigonality algorithm}

We describe in more detail the algorithm to detect and compute the trigonality of a curve. In particular, we explain now the computation of a threefold map from $C$ to $\pj^1$. Let $S$ be the surface intersection of quadrics containing $\varphi(C)$ and $S\rightarrow\pj^2$ the parametrization obtained with the Lie algebra method. In all cases we will obtain a threefold map from $\varphi(C)$ to $\pj^1$ which can be pulled back to $C$.

\begin{itemize}
\item If $g=3$, the canonical curve is a smooth quartic in $\pj^2$. In this case, one can easily compute a $g^1_3$ on $\varphi(C)$, at least in theory: it suffices to take a point $p$ on it and consider the pencil of lines through it, since each one intersects $\varphi(C)$ in $p$ and three more points. In practice, finding a point with coefficients in the base field is problematic, unless one accepts working on algebraic extensions.

\item If $S$ is a rational normal scroll not isomorphic to $\pj^1\times\pj^1$, we can compute the map in \eqref{EQ: isomorphism to P1} explicitly.

\item If $S$ is isomorphic to $\pj^1\times\pj^1$, we compute a map $\rho:\pj^1\times\pj^1\rightarrow\pj^1$ and form the composition
\[C\stackrel{\varphi}{\rightarrow}\varphi(C)\stackrel{i}{\hookrightarrow}S\stackrel{\cong}{\rightarrow}\pj^1\times\pj^1\stackrel{\rho}{\rightarrow}\pj^1.\]

\item If $C$ is a plane quintic, it is not trigonal.
\end{itemize}

\section{Computational experiences}

We have tested our Magma V2.14-7 \cite{MAGMA} implementation against many examples of trigonal curves. The computer used is a 64 Bit, Dual AMD Opteron Processor 250 (2.4 GHZ) with 8 GB RAM. We have generated trigonal curves in the following two ways:

\begin{enumerate}
 \item Let $C:f(x,y,z)=0$ with $\deg_y f=3$. Then the projection $(x:y:z)\mapsto(x:z)$ is a $3:1$ map to $\pj^1$. The genus of a polynomial of degree 3 in $y$ and degree $d$ in $x$ is $2(d-1)$ generically. The size of the coefficients is controlled directly.
 \item Let $C$ be defined by the affine equation $\Resultant_u(F,G)=0$ where
\[\begin{array}{l@{\ }l}
0 = x^3-a_1(u)x-a_2(u) & =:F, \\
0 = y-a_3(u)-a_4(u)x-a_5(u)x^2 & =:G
\end{array}\]
for some polynomials $a_1,\ldots,a_5$. This clearly gives a field extension of degree 3, thus there is a $3:1$ map from $C$ to the affine line. The degree and coefficient size for a given genus are significantly larger than for the previous construction.
\end{enumerate}

These are our timed results\footnoteremember{myfootnote}{Last minute improvements in our implementation have reduced the running times by a factor of about 5, for many of the entries in the tables.} for samples of ten random polynomials of different degrees, genera and coefficient sizes.

\[\begin{array}{cc|ccc}
\deg_x & \mathrm{bit\ height} & \mathrm{genus} & \deg & \mathrm{seconds} \\ \hline
3 & 5 & 4 & 6 & 0.5-0.65 \\
3 & 50 & 4 & 6 & 2.09-2.27 \\
6 & 5 & 10 & 9 & 14-17 \\
6 & 50 & 10 & 9 & 54-61 \\
10 & 5 & 18 & 13 & 271-342 \\
10 & 50 & 18 & 13 & 1059-1193 \\
15 & 5 & 28 & 18 & 3477-5317 \\
\end{array}\]

For the second method, we choose $a_1,\ldots,a_5$ randomly of degree $d$ and maximum coefficient size $e$. These are the time results\footnoterecall{myfootnote} for samples of ten random polynomials, for different values of $d,e$.

\[\begin{array}{c|cccc}
(d,e) & \mathrm{genus} & \deg & \mathrm{bit\ height} & \mathrm{seconds} \\ \hline
(4,2) & 4 & 15-20 & 9-16 & 18-62 \\
(4,10) & 4 & 20 & 26-35 & 87-191 \\
(5,2) & 4-6 & 20-25 & 17-21 & 162-2353  \\
(5,10) & 4-6 & 23-25 & 34-41 & 1334-7940 \\
(6,2) & 6-7 & 25-30 & 20-24 & 2992-22650 \\
\end{array}\]

Our Magma implementation can be obtained by contacting us directly: \texttt{josef.schicho@oeaw.ac.at} and \texttt{david.sevilla@oeaw.ac.at}.

\bibliographystyle{plain}

\begin{thebibliography}{9}

\bibitem{Enriques1919}
F. Enriques, \emph{Sulle curve canoniche di genere p dello spazio a p-1 dimensioni}. Rend. dell'Acc. delle Scienze di Bologna 23 (1919), 80--82.

\bibitem{Babbage1939}
D.W. Babbage, \emph{A note on the quadrics through a canonical curve}. J. London Math Soc. 14 (1939), 310--315.

\bibitem{Petri1923}
K. Petri, \emph{\"Uber die invariante Darstellung algebraischer Funktionen einer Ver\"anderlichen [in German]}. Math. Ann. 88 (1923), no. 3-4, 242--289.

\bibitem{GriffithsHarris1978}
P. Griffiths and J. Harris, \emph{Principles of algebraic geometry}.
Pure and Applied Mathematics. Wiley-Interscience [John Wiley \& Sons], New York, 1978. ISBN: 0-471-32792-1.

\bibitem{Schicho2008}
W. A. de Graaf, J. P\'ilnikov\'a and J. Schicho, \emph{Parametrizing Del Pezzo surfaces of degree 8 using Lie algebras}. J. Symb. Comp. 44 (2009) no. 1, 1--14.

\bibitem{Schicho2006}
W. A. de Graaf, M. Harrison, J. P{\'i}nikov{\'a} and J. Schicho, \emph{A Lie algebra method for the parametrization of Severi-Brauer surfaces}. J. Algebra 303 (2006) no.2, 514--529.

\bibitem{MAGMA}
W. Bosma, J. Cannon and C. Playoust, \emph{The Magma algebra system. I. The user language}. J. Symbolic Comput. 24 (1997), no. 3-4, 235--265.

\end{thebibliography}

\end{document}